\documentclass[eqthmnum]{jt-calcs}
\usepackage{tikz}

\usepackage[backend=bibtex,style=alphabetic,sorting=nyt,minnames=2]{biblatex}
\bibliography{bibliography}
\AtEveryBibitem{%
  \clearfield{pagetotal}%
}

\title{Induced Characters of Type $\D$ Weyl Groups and the Littlewood--Richardson Rule}
\author{Jay Taylor}
\address{FB Mathematik, TU Kaiserslautern, Postfach 3049, 67653 Kaiserslautern, Germany}
\email{taylor@mathematik.uni-kl.de}
\mscno{2010}{20C15}{20G40}
\keywords{Weyl groups, type D, reflection subgroups, induced characters}

\begin{document}
\begin{abstract}
For any ordinary irreducible character of a maximal reflection subgroup of type $\D_a\D_b$ of a type $\D$ Weyl group we give an explicit decomposition of the induced character in terms of Littlewood--Richardson coefficients.
\end{abstract}

\section{Introduction}
\begin{pa}
Let $W$ be a finite Weyl group and $H \leqslant W$ a reflection subgroup. A natural question that occurs in the representation theory of $W$ is to describe the decomposition of the induced character $\Ind_H^W(\chi)$ for any (ordinary) irreducible character $\chi$ of $H$. By a general process one can reduce this to the case where $W$ is irreducible and $H$ is maximal. The classical situation occurs when $W$ is the symmetric group $\mathfrak{S}_n$ then $H$ is a Young subgroup $\mathfrak{S}_a\mathfrak{S}_b$ with $a+b=n$. In this case the decomposition is given by the classical Littlewood--Richardson coefficients which are combinatorially computable using the Littlewood--Richardson rule. For exceptional type Weyl groups one can simply tackle this using the CHEVIE package \cite{geck-hiss:1996:CHEVIE}.
\end{pa}

\begin{pa}
The remaining cases occur when $W$ is of type $\B_n$ or $\D_n$. In the first case we have $H$ is either of type $\D_n$ or $\B_a\B_b$ with $a+b= n$ (see \cite[Lemma 5]{tumarkin-felikson:2005:reflection-subgroups}). Both of these cases are treated, for instance, in the book of Geck and Pfeiffer (see \cite[5.6.1, 6.1.3]{geck-pfeiffer:2000:characters-of-finite-coxeter-groups}). Now if $W$ is of type $\D_n$ then we have $H$ is either of type $\A_{n-1}$, $\D_{n-1}$ or $\D_a\D_b$ with $2\leqslant a,b \leqslant n-2$ and $a+b=n$ (see \cite[Corollary 1]{tumarkin-felikson:2005:reflection-subgroups}). For the so-called non-degenerate irreducible characters of $W$ the multiplicities are easily deduced from the type $\B$ case. However the multiplicities for degenerate characters requires some work. It is the purpose of this note to consider the case where $H$ is of type $\D_{n-1}$ or $\D_a\D_b$. Here we are able to get a complete description of the multiplicities in terms of Littlewood--Richardson coefficients. We suspect that this result is well known to the experts but we failed to find it in the literature. When $H$ is of type $\A_{n-1}$ then certain special cases are considered in \cite[5.6.3]{geck-pfeiffer:2000:characters-of-finite-coxeter-groups} and \cite[3.7]{geck:2015:on-kottwitz-conjecture} but a combinatorial formula in the general case is still to be obtained.
\end{pa}

\begin{pa}
We briefly point out why such a result is relevant in a wider context and also the motivation for the author. Let $\bG$ be a connected reductive algebraic group over an algebraically closed field of characteristic $p$ and $G = \bG^F$ the fixed points under a Frobenius endomorphism $F : \bG \to \bG$. In the representation theory of $G$ the representations of Weyl groups occur frequently and in many varied ways. In \cite{lusztig:1986:on-the-character-values} Lusztig showed that the restriction of a character sheaf of $\bG$ to the unipotent variety could be explicitly decomposed in terms of certain intersection cohomology complexes, which are defined on the closure of a unipotent conjugacy class. The coefficients in this decomposition correspond precisely to the coefficients in the decomposition of an induced character from a subgroup of a relative Weyl group of $\bG$, which is a reflection subgroup when $Z(\bG)$ is connected. Recently the author has extended the work of Lusztig to obtain an explicit formula for the value of the characteristic function of a character sheaf at a unipotent element (see \cite{taylor:2014:evaluating-characteristic-functions}). In this formula such multiplicities also occur, as well as their $F$-coset analogues.
\end{pa}

\begin{pa}
If $\bG$ is an adjoint simple group of type $\B_n$ then under certain mild restrictions Lusztig was able to prove, for unipotently supported character sheaves, his conjecture relating characteristic functions of character sheaves to the irreducible characters of $G$ (see \cite[5.3(c)]{lusztig:1986:on-the-character-values} and \cite[13.7]{lusztig:1984:characters-of-reductive-groups}). A key role in the proof of Lusztig's result was the type $\B$ analogue of \cref{lem:typeD-GP-Lemm6.1.3}, as given in \cite[6.1.3]{geck-pfeiffer:2000:characters-of-finite-coxeter-groups}. It is the authors hope to extend Lusztig's result to the case where $\bG$ is of type $\D_n$, where one needs the information contained in \cref{lem:typeD-GP-Lemm6.1.3}.
\end{pa}

\begin{acknowledgments}
The author gratefully acknowledges the financial support of ERC Advanced Grant 291512 awarded to Prof.\ Gunter Malle. He also warmly thanks the referee for their careful reading of the manuscript and useful suggestions.
\end{acknowledgments}

\section{The Result}
\begin{pa}[Notation]\label{pa:notation}
Let $G$ be a finite Weyl group. We assume that $\mathbb{Q} \subseteq \mathbb{K}$ is a fixed field extension, then the group algebra $\mathbb{K}G$ is semisimple; we denote by $\Irr(G)$ the irreducible characters of $G$ over $\mathbb{K}$. If $\chi, \psi \in \Irr(G)$ are irreducible characters then we set
\begin{equation*}
\langle \chi,\psi\rangle_G = \dim\Hom_{\mathbb{K}G}(M_{\chi},M_{\psi}),
\end{equation*}
where $M_{\chi}$ (resp.\ $M_{\psi}$) is a $\mathbb{K}G$-module affording $\chi$ (resp.\ $\psi$). We extend this linearly to a bilinear map $\langle -,-\rangle_G : \mathbb{Z}\Irr(G) \times \mathbb{Z}\Irr(G) \to \mathbb{Z}$. Given any element $g \in G$ we will denote by $\Cl_G(g)$ (resp.\ $C_G(g)$) the conjugacy class (resp.\ centraliser) of $g$. If $G$ has reflection representation $G \to \GL(V)$ then, for any subgroup $H \leqslant G$, we denote by $j_H^G$ the Lusztig--Macdonald--Spaltenstein induction with respect to $V$ as defined in \cite[\S5.2]{geck-pfeiffer:2000:characters-of-finite-coxeter-groups}.

By a \emph{partition} we will mean a finite (possibly empty) weakly decreasing sequence of strictly positive integers $\alpha = (\alpha_1,\dots,\alpha_k)$. For any such partition $\alpha$ we set
\begin{equation*}
|\alpha| = \begin{cases}
0 &\text{if }\alpha = ()\\
\alpha_1+\cdots+\alpha_k &\text{if }\alpha\neq().
\end{cases}
\end{equation*}
We say $\alpha$ is a \emph{partition of $n\geqslant 0$} if $|\alpha| = n$. By a \emph{bipartition of $n$} we will refer to an ordered pair of partitions $(\alpha;\beta)$ such that $|\alpha|+|\beta| = n$. For any $n \geqslant 0$ we will denote by $\mathcal{P}(n)$ (resp.\ $\mathcal{BP}(n)$) the set of partitions (resp.\ bipartitions) of $n$. We will also denote by $\mathcal{X}(n)$ the set of all pairs of integers $(a,b)$ such that $a,b > 0$ and $a+b=n$. If $n$ is not a positive integer then we set $\mathcal{P}(n) = \mathcal{BP}(n) = \mathcal{X}(n) = \emptyset$. Given any partition $\alpha = (\alpha_1,\dots,\alpha_k)$ we write $\ell(\alpha) = k$ for the length of the partition, which is simply 0 if $\alpha$ is empty. If $\beta = (\beta_1,\dots,\beta_{\ell})$ is a second partition then we write $\alpha\cup\beta \in \mathcal{P}(|\alpha|+|\beta|)$ for the partition obtained by reordering the sequence $(\alpha_1,\dots,\alpha_k,\beta_1,\dots,\beta_{\ell})$. Note we clearly have $\ell(\alpha\cup\beta) = \ell(\alpha)+\ell(\beta)$.
\end{pa}

\begin{pa}
For any $n \geqslant 0$ let $\widetilde{W}_n$ be a Coxeter group of type $\B_n$ with generators $\{s_1,\dots,s_{n-1},t_n\}$ and corresponding Coxeter diagram
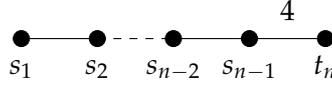
\begin{center}
\begin{tikzpicture}[text height=1.5ex,text depth=.25ex]
\clip (-.5,-.75) rectangle (4.5,.5);

\filldraw (0,0) circle [radius=3pt];
\filldraw (1,0) circle [radius=3pt];
\filldraw (2,0) circle [radius=3pt];
\filldraw (3,0) circle [radius=3pt];
\filldraw (4,0) circle [radius=3pt];
\draw (0,0) -- (1,0);
\draw[dashed] (1,0) -- (2,0);
\draw (2,0) -- (3,0);
\draw (3,0) -- (4,0) node [midway,above=2pt] {4};

\node at (0,-.35) {$s_1$};
\node at (1,-.35) {$s_2$};
\node at (2,-.35) {$s_{n-2}$};
\node at (3,-.35) {$s_{n-1}$};
\node at (4,-.35) {$t_n$};
\end{tikzpicture}
\end{center}
For $0 \leqslant n \leqslant 1$ we set $W_n \leqslant \widetilde{W}_n$ to be the trivial subgroup. Assume now that $n > 1$ then setting $s_n = t_ns_{n-1}t_n$ we have the subgroup $W_n = \langle s_1,\dots,s_n\rangle \leqslant \widetilde{W}_n$ is naturally a Coxeter group of type $\D_n$ with corresponding Coxeter diagram
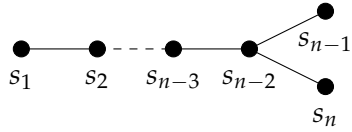
\begin{center}
\begin{tikzpicture}[text height=1.5ex,text depth=.25ex]
\clip (-.5,-1.2) rectangle (4.5,1);

\filldraw (0,0) circle [radius=3pt];
\filldraw (1,0) circle [radius=3pt];
\filldraw (2,0) circle [radius=3pt];
\filldraw (3,0) circle [radius=3pt];
\filldraw (4,.5) circle [radius=3pt];
\filldraw (4,-.5) circle [radius=3pt];
\draw (0,0) -- (1,0);
\draw[dashed] (1,0) -- (2,0);
\draw (2,0) -- (3,0);
\draw (3,0) -- (4,.5);
\draw (3,0) -- (4,-.5);

\node at (0,-.35) {$s_1$};
\node at (1,-.35) {$s_2$};
\node at (2,-.35) {$s_{n-3}$};
\node at (3,-.35) {$s_{n-2}$};
\node at (4,.15) {$s_{n-1}$};
\node at (4,-.85) {$s_n$};
\end{tikzpicture}
\end{center}
Note that $W_2$ is of type $\D_2 = \A_1\A_1$ and $W_3$ is of type $\D_3 = \A_3$.
\end{pa}

\begin{pa}
We assume now that $n \geqslant 1$. Recall that there is a bijection $\mathcal{BP}(n) \to \Irr(\widetilde{W}_n)$ which we denote by $(\alpha;\beta) \mapsto [\alpha;\beta]$. We fix this bijection to be the one defined in \cite[5.5.6]{geck-pfeiffer:2000:characters-of-finite-coxeter-groups}. The group $W_n \leqslant \widetilde{W}_n$ is an index 2 normal subgroup. In particular, the restriction $\overline{[\alpha;\beta]} = \Res_{W_n}^{\widetilde{W}_n}([\alpha;\beta])$ of any irreducible character $[\alpha;\beta] \in \Irr(\widetilde{W}_n)$ to $W_n$ is either irreducible or the sum of two irreducible characters. From \cite[5.6.1]{geck-pfeiffer:2000:characters-of-finite-coxeter-groups} we recall that
\begin{equation*}
\overline{[\alpha;\beta]} = \overline{[\beta;\alpha]}
\end{equation*}
is an irreducible character of $W_n$ if and only if $\alpha \neq \beta$. Assuming $\alpha \neq \beta$ then we call the irreducible character $\overline{[\alpha;\beta]} \in \Irr(W_n)$ \emph{non-degenerate}. Now assume $n$ is even and $\alpha \in \mathcal{P}(n/2)$ is a partition. Then we have
\begin{equation*}
\overline{[\alpha;\alpha]} = \overline{[\alpha;\alpha]}_+ + \overline{[\alpha;\alpha]}_-,
\end{equation*}
where $\overline{[\alpha;\alpha]}_{\pm} \in \Irr(W_n)$. We say $\overline{[\alpha;\alpha]}_{\pm}$ is a \emph{degenerate} character of $W_n$ of type $\pm$. It will be convenient to also write $\overline{[\alpha;\beta]}_{\pm}$ even when $\alpha \neq \beta$; in this case the sign should simply be ignored. To conclude we have the irreducible characters of $W_n$ are
\begin{equation*}
\Irr(W_n) = \{\overline{[\alpha;\beta]} \mid (\alpha;\beta) \in \mathcal{BP}(n)\text{ and }\alpha\neq\beta\}\cup\{\overline{[\alpha;\alpha]}_{\pm} \mid \alpha \in \mathcal{P}(n/2)\}.
\end{equation*}
For any $X \in \Irr(W_n)$ we set
\begin{equation*}
\varepsilon(X) = \begin{cases}
0 &\text{if }X\text{ is non-degenerate}\\
\pm 1 &\text{if }X\text{ is degenerate of type }\pm.
\end{cases}
\end{equation*}
\end{pa}

\begin{pa}\label{pa:iso-symmetric}
To make the above parameterisation of $\Irr(W_n)$ concrete we must distinguish the irreducible characters $\overline{[\alpha;\alpha]}_{\pm}$ which we do following \cite[5.6.3]{geck-pfeiffer:2000:characters-of-finite-coxeter-groups}. Let
\begin{equation*}
H_n^+ = \langle s_1,\dots,s_{n-2},s_{n-1} \rangle \qquad\qquad H_n^- = \langle s_1,\dots,s_{n-2},s_n\rangle
\end{equation*}
then both of these are subgroups of $W_n$ isomorphic to the symmetric group $\mathfrak{S}_n$ on $\{1,\dots,n\}$. If $n$ is even then they are not conjugate. Recall that we have a bijection $\mathcal{P}(n) \to \Irr(\mathfrak{S}_n)$ defined as in \cite[5.4.7]{geck-pfeiffer:2000:characters-of-finite-coxeter-groups}; we denote this by $\alpha \mapsto [\alpha]$. We have a natural isomorphism $H_n^{\pm} \cong \mathfrak{S}_n$ given by $s_i \mapsto (i,i+1)$, for $1 \leqslant i \leqslant n-2$, and $s_{n-1}, s_n \mapsto (n-1,n)$. Under this isomorphism we will identify the set of characters $\Irr(H_n^{\pm})$ with $\Irr(\mathfrak{S}_n)$. Now for any $\alpha \in \mathcal{P}(n/2)$ we set
\begin{equation*}
\overline{[\alpha;\alpha]}_+ = j_{H_n^+}^{W_n}([\alpha\cup\alpha]) \qquad\qquad \overline{[\alpha;\alpha]}_- = j_{H_n^-}^{W_n}([\alpha\cup\alpha])
\end{equation*}
and this distinguishes the characters (c.f.\ \cref{pa:notation}).
\end{pa}

\begin{pa}\label{pa:reflection-subgroups-WaWb}
Assume $n \geqslant 4$ then setting $s_0 = t_0s_1t_0 \in W_n$ with
\begin{equation*}
t_0 = s_1s_2\cdots s_{n-1}t_ns_{n-1}\cdots s_2s_1
\end{equation*}
we obtain the extended Coxeter diagrams
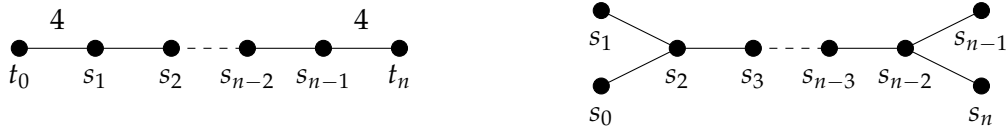
\begin{center}
\begin{tikzpicture}[text height=1.5ex,text depth=.25ex]
\clip (-1.5,-1.1) rectangle (4.5,1);

\filldraw (-1,0) circle [radius=3pt];
\filldraw (0,0) circle [radius=3pt];
\filldraw (1,0) circle [radius=3pt];
\filldraw (2,0) circle [radius=3pt];
\filldraw (3,0) circle [radius=3pt];
\filldraw (4,0) circle [radius=3pt];
\draw (-1,0) -- (0,0) node [midway,above=2pt] {4};
\draw (0,0) -- (1,0);
\draw[dashed] (1,0) -- (2,0);
\draw (2,0) -- (3,0);
\draw (3,0) -- (4,0) node [midway,above=2pt] {4};

\node at (-1,-.35) {$t_0$};
\node at (0,-.35) {$s_1$};
\node at (1,-.35) {$s_2$};
\node at (2,-.35) {$s_{n-2}$};
\node at (3,-.35) {$s_{n-1}$};
\node at (4,-.35) {$t_n$};
\end{tikzpicture}
\qquad\qquad
\begin{tikzpicture}[text height=1.5ex,text depth=.25ex]
\clip (-1.5,-1.1) rectangle (4.5,1);

\filldraw (-1,.5) circle [radius=3pt];
\filldraw (-1,-.5) circle [radius=3pt];
\filldraw (0,0) circle [radius=3pt];
\filldraw (1,0) circle [radius=3pt];
\filldraw (2,0) circle [radius=3pt];
\filldraw (3,0) circle [radius=3pt];
\filldraw (4,.5) circle [radius=3pt];
\filldraw (4,-.5) circle [radius=3pt];

\draw (-1,.5) -- (0,0);
\draw (-1,-.5) -- (0,0);
\draw (0,0) -- (1,0);
\draw[dashed] (1,0) -- (2,0);
\draw (2,0) -- (3,0);
\draw (3,0) -- (4,.5);
\draw (3,0) -- (4,-.5);

\node at (-1,-.85) {$s_0$};
\node at (-1,.15) {$s_1$};
\node at (0,-.35) {$s_2$};
\node at (1,-.35) {$s_3$};
\node at (2,-.35) {$s_{n-3}$};
\node at (3,-.35) {$s_{n-2}$};
\node at (4,.15) {$s_{n-1}$};
\node at (4,-.85) {$s_n$};
\end{tikzpicture}
\end{center}
We set $\widetilde{\mathbb{S}}_0 = \{t_0,s_1,\dots,s_{n-1},t_n\}$ and $\mathbb{S}_0 = \{s_0,\dots,s_n\}$. Given $(a,b) \in \mathcal{X}(n)$ we have corresponding reflection subgroups $\widetilde{W}_a\widetilde{W}_b = \langle\widetilde{\mathbb{S}}_0\setminus\{s_a\}\rangle$ of $\widetilde{W}_n$ and
\begin{equation*}
W_aW_b = \begin{cases}
\langle\mathbb{S}_0\setminus\{s_0,s_1\}\rangle &\text{if }a = 1\\
\langle\mathbb{S}_0\setminus\{s_{n-1},s_n\}\rangle &\text{if }a = n-1\\
\langle\mathbb{S}_0\setminus\{s_a\}\rangle &\text{otherwise}
\end{cases}
\end{equation*}
of $W_n$. Note that we clearly have $W_aW_b = \widetilde{W}_a\widetilde{W}_b \cap W_n$. Setting $H_a^{\pm} = W_a \cap H_n^{\pm}$ determines two copies of $\mathfrak{S}_a$ inside $W_a$ (similarly for $W_b$). With respect to these choices the degenerate characters of $W_a$ and $W_b$ are concretely labelled as in \cref{pa:iso-symmetric}.
\end{pa}

\begin{pa}
Assume $(a,b) \in \mathcal{X}(n)$ and $\alpha \in \mathcal{P}(a)$, $\beta \in \mathcal{P}(b)$ and $\gamma \in \mathcal{P}(n)$ are partitions. We define the corresponding Littlewood--Richardson coefficient to be
\begin{equation*}
c_{\alpha\beta}^{\gamma} = \langle \Ind_{\mathfrak{S}_a\mathfrak{S}_b}^{\mathfrak{S}_n}([\alpha]\boxtimes[\beta]), [\gamma] \rangle_{\mathfrak{S}_n}.
\end{equation*}
Here $\mathfrak{S}_a\mathfrak{S}_b$ denotes the Young subgroup of $\mathfrak{S}_n$ preserving the sets $\{1,\dots,a\}$ and $\{a+1,\dots,n\}$. Furthermore, $[\alpha]\boxtimes[\beta]$ denotes the tensor product of the characters $[\alpha]$ and $[\beta]$ so that $([\alpha]\boxtimes[\beta])(wx) = [\alpha](w)[\beta](x)$ for all $wx \in \mathfrak{S}_a\mathfrak{S}_b$. Now let $\alpha = (\alpha_1;\alpha_2) \in \mathcal{BP}(a)$, $\beta = (\beta_1;\beta_2) \in \mathcal{BP}(b)$ and $\gamma = (\gamma_1;\gamma_2) \in \mathcal{BP}(n)$ be bipartitions. We will denote by $\mathfrak{S}_{\alpha}$ the symmetric group $\mathfrak{S}_1$ (resp.\ $\mathfrak{S}_2$) if $\alpha_1 = \alpha_2$ (resp.\ $\alpha_1\neq\alpha_2$) and similarly for $\beta$. With this we define
\begin{equation*}
a_{\alpha\beta}^{\gamma} = \sum_{\sigma \in \mathfrak{S}_{\alpha}}\sum_{\tau \in \mathfrak{S}_{\beta}} c_{\alpha_{\sigma(1)}\beta_{\tau(1)}}^{\gamma_1}c_{\alpha_{\sigma(2)}\beta_{\tau(2)}}^{\gamma_2}
\end{equation*}
where the sum is over all permutations such that $|\gamma_i| = |\alpha_{\sigma(i)}| + |\beta_{\tau(i)}|$ for $i \in \{1,2\}$. If no such permutations exist then we have $a_{\alpha\beta}^{\gamma} = 0$. Note that this definition does not depend upon the ordering of the partitions in the bipartitions $\alpha$, $\beta$ or $\gamma$. With this notation in place we may now state our result.
\end{pa}

\begin{prop}\label{lem:typeD-GP-Lemm6.1.3}
Assume $n \geqslant 4$ and $(a,b) \in \mathcal{X}(n)$. Let $E = A \boxtimes B \in \Irr(W_aW_b)$ be an irreducible character then we denote by $\alpha = (\alpha_1;\alpha_2) \in \mathcal{BP}(a)$ and $\beta = (\beta_1;\beta_2) \in \mathcal{BP}(b)$ bipartitions such that $A = \overline{[\alpha_1;\alpha_2]}_{\pm}$ and $B = \overline{[\beta_1;\beta_2]}_{\pm}$. Now let $X \in \Irr(W_n)$ be an irreducible character and let $\gamma = (\gamma_1;\gamma_2) \in \mathcal{BP}(n)$ be a bipartition such that $X = \overline{[\gamma_1;\gamma_2]}_{\pm}$. Then
\begin{equation*}
\langle\Ind_{W_aW_b}^{W_n}(E),X\rangle_{W_n} = \begin{cases}
a_{\alpha\beta}^{\gamma} &\text{if }X\text{ is non-degenerate},\\
\frac{1}{2}(a_{\alpha\beta}^{\gamma} + \varepsilon(A)\varepsilon(B)\varepsilon(X)c_{\alpha_1\beta_1}^{\gamma_1}) &\text{if }X\text{ is degenerate}.
\end{cases}
\end{equation*}
\end{prop}

\begin{rem}
We note that in the interesting case when $\varepsilon(A)\varepsilon(B)\varepsilon(X) \neq 0$, i.e.\ when $\alpha_1=\alpha_2$, $\beta_1 = \beta_2$ and $\gamma_1 = \gamma_2$, then the result may be written as
\begin{equation*}
\langle\Ind_{W_aW_b}^{W_n}(E),X\rangle_{W_n} = \frac{1}{2}c_{\alpha_1\beta_1}^{\gamma_1}(c_{\alpha_1\beta_1}^{\gamma_1} + \varepsilon(A)\varepsilon(B)\varepsilon(X))
\end{equation*}
because $a_{\alpha\beta}^{\gamma} = (c_{\alpha_1\beta_1}^{\gamma_1})^2$.
\end{rem}

\begin{proof}
For notational convenience we set $H = W_aW_b$, $\widetilde{H} = \widetilde{W}_a\widetilde{W}_b$, $W = W_n$ and $\widetilde{W} = \widetilde{W}_n$. By the Clifford theoretic description of the irreducible characters of $W$ and \cite[6.1.3]{geck-pfeiffer:2000:characters-of-finite-coxeter-groups} it is easy to deduce that
\begin{equation}\label{eq:type-B-mult}
\langle \Ind_{H}^{\widetilde{W}}(E), [\gamma_1;\gamma_2]\rangle_{\widetilde{W}} = \langle \Ind_{\widetilde{H}}^{\widetilde{W}}\Ind_{H}^{\widetilde{H}}(E), [\gamma_1;\gamma_2]\rangle_{\widetilde{W}} = a_{\alpha\beta}^{\gamma}.
\end{equation}
Applying Frobenius reciprocity we obtain the multiplicity of any non-degen\-erate character in $\Ind_H^W(E)$. Thus we are left with considering the case of degenerate characters.

Assume now that $X \in \Irr(W)$ is degenerate so $n$ is even. If $B$ is non-degenerate then we have
\begin{equation*}
{}^{s_n}\Ind_{H}^{W}(E) = \Ind_{H}^{W}({}^{s_n}E) = \Ind_{H}^{W}(E).
\end{equation*}
As conjugation by $s_n$ permutes the degenerate characters $\overline{[\gamma_1;\gamma_2]}_{\pm}$ we get from \cref{eq:type-B-mult} that
\begin{equation*}
\langle \Ind_{H}^{W}(E),X\rangle_{W} = \langle \Ind_{H}^{W}(E),{}^{s_n}X\rangle_{W} = \frac{1}{2}a_{\alpha\beta}^{\gamma}.
\end{equation*}
Similarly, if $A$ is non-degenerate we may conjugate the induced character by the element $s_0$ to deduce the same result. Hence we need only deal with the situation where both $A$ and $B$ are degenerate, which we now assume to be the case. In particular both $a$ and $b$ are even.

We denote by $\Pi_X$ (resp.\ $\Delta_X$) the character $\overline{[\gamma_1;\gamma_2]} = \overline{[\gamma_1;\gamma_2]}_++\overline{[\gamma_1;\gamma_2]}_-$ (resp.\ the difference character $\overline{[\gamma_1;\gamma_2]}_+ - \overline{[\gamma_1;\gamma_2]}_-$). Using the fact that any degenerate character $X \in \Irr(W)$ may be written as $\frac{1}{2}(\Pi_X+\varepsilon(X)\Delta_X)$ we have
\begin{equation*}
\langle \Ind_{H}^{W}(E), X \rangle_{W} = \frac{1}{2}\langle \Ind_{H}^{W}(E),\Pi_X\rangle_{W} + \frac{1}{2}\varepsilon(X)\langle \Ind_{H}^{W}(E),\Delta_X\rangle_{W}.
\end{equation*}
Note that $\langle \Ind_{H}^{W}(E),\Pi_X\rangle_{W} = a_{\alpha\beta}^{\gamma}$ by \cref{eq:type-B-mult} so we have only to compute the multiplicity $\langle \Ind_{H}^{W}(E),\Delta_X\rangle_{W}$. As $A$ and $B$ are both degenerate we can apply this approach to the character $E$ to obtain
\begin{align*}
\Ind_{H}^{W}(E) &= \frac{1}{4}\Big[ \Ind_{H}^{W}(\Pi_A\boxtimes\Pi_B) + \varepsilon(B)\Ind_{H}^{W}(\Pi_A\boxtimes\Delta_B)\\
&\qquad\qquad+ \varepsilon(A)\Ind_{H}^{W}(\Delta_A\boxtimes\Pi_B) + \varepsilon(A)\varepsilon(B)\Ind_{H}^{W}(\Delta_A\boxtimes\Delta_B)\Big].
\end{align*}
We now observe that the character $\Pi_A$ (resp.\ $\Pi_B$) is invariant under conjugation by $s_0$ (resp.\ $s_n$). Hence, applying the previous argument we see that the degenerate characters $\overline{[\gamma_1;\gamma_2]}_{\pm}$ occur with the same multiplicity in the first three induced characters. This implies that their inner product with the difference character $\Delta_X$ is 0. In particular
\begin{equation*}
\langle \Ind_{H}^{W}(E),\Delta_X\rangle_{W} = \frac{1}{4}\varepsilon(A)\varepsilon(B)\langle \Ind_{H}^{W}(\Delta_A\boxtimes\Delta_B),\Delta_X\rangle_{W}.
\end{equation*}
For notational convenience we will set $\Theta = \Ind_{H}^{W}(\Delta_A\boxtimes\Delta_B)$.

We now wish to compute the multiplicity on the right hand side but to do this we will need some information about conjugacy classes and the difference character $\Delta_X$. For any partition $\pi = (\pi_1,\dots,\pi_k) \in \mathcal{P}(n/2)$ we will denote by $w_{2\pi}^+ \in H_n^+ \leqslant W$ an element of cycle type $2\pi = (2\pi_1,\dots,2\pi_k)$. This definition makes sense because of the isomorphism chosen in \cref{pa:iso-symmetric}. We also denote by $w_{2\pi}^-$ the conjugate $s_nw_{2\pi}^+s_n \in H_n^- \leqslant W$ and by $w_{\pi} \in \mathfrak{S}_{n/2}$ an element of cycle type $\pi$.

Now, according to \cite[10.4.9]{geck-pfeiffer:2000:characters-of-finite-coxeter-groups} and \cite[3.5]{geck:2015:on-kottwitz-conjecture}, we have for all $w \in W$ that
\begin{equation*}
\Delta_X(w) = \begin{cases}
\pm (-1)^{n/2}2^{\ell(\pi)}[\gamma_1](w_{\pi}) &\text{if }w \in \Cl_{W}(w_{2\pi}^{\pm})\text{ for some }\pi\in\mathcal{P}(n/2)\\
0 &\text{otherwise}.
\end{cases}
\end{equation*}
Here $[\gamma_1]$ denotes an irreducible character of the symmetric group $\mathfrak{S}_{n/2}$. Note that the statement in \cite[10.4.10]{geck-pfeiffer:2000:characters-of-finite-coxeter-groups} identifying the characters in \cite[10.4.6]{geck-pfeiffer:2000:characters-of-finite-coxeter-groups} with those in \cite[5.6.3]{geck-pfeiffer:2000:characters-of-finite-coxeter-groups} is not correct; the correct statement is given in \cite[3.5]{geck:2015:on-kottwitz-conjecture}. We will also need the following identities concerning centraliser orders which are easily obtained by direct computation
\begin{equation*}
|C_{W}(w_{2\pi}^{\pm})| = |C_{\widetilde{W}}(w_{2\pi}^{\pm})| = 2^{\ell(\pi)}|C_{H_n^{\pm}}(w_{2\pi}^{\pm})| = 2^{2\ell(\pi)}|C_{\mathfrak{S}_{n/2}}(w_{\pi})|.
\end{equation*}

To compute the multiplicity $\langle \Theta,\Delta_X\rangle_{W}$ we will need to determine the value of the induced character at $w_{2\pi}^{\pm}$, for which we use the explicit induction formula given in \cite[pg.\ 64]{isaacs:2006:character-theory-of-finite-groups}. For this we need to determine the orbits of $H$ acting on $H \cap \Cl_{W}(w_{2\pi}^{\pm})$ by conjugation. With this in mind we define
\begin{equation*}
\mathcal{B}_{\pi} = \{(\delta,\epsilon) \in \mathcal{P}(a/2) \times \mathcal{P}(b/2) \mid \delta\cup\epsilon = \pi\},
\end{equation*}
(c.f.\ \cref{pa:notation}). We then have
\begin{equation*}
\begin{cases}
\{w_{2\delta}^+w_{2\epsilon}^+,w_{2\delta}^-w_{2\epsilon}^- \mid (\delta,\epsilon) \in \mathcal{B}_{\pi}\} &\text{for $w_{2\pi}^+$}\\
\{w_{2\delta}^+w_{2\epsilon}^-,w_{2\delta}^-w_{2\epsilon}^+ \mid (\delta,\epsilon) \in \mathcal{B}_{\pi}\} &\text{for $w_{2\pi}^-$}
\end{cases}
\end{equation*}
is a complete set of representatives for the desired orbits where
\begin{equation*}
w_{2\delta}^- = s_0w_{2\delta}^+s_0 \in H_a^- \qquad\text{and}\qquad w_{2\epsilon}^- = s_nw_{2\epsilon}^+s_n \in H_b^-.
\end{equation*}
In addition, $\{w_{\delta}w_{\epsilon} \mid (\delta,\epsilon) \in \mathcal{B}_{\pi}\} \subseteq \mathfrak{S}_{a/2}\mathfrak{S}_{b/2}$ is a complete set of representatives for the orbits of $\mathfrak{S}_{a/2}\mathfrak{S}_{b/2}$ acting on $\mathfrak{S}_{a/2}\mathfrak{S}_{b/2}\cap\Cl_{\mathfrak{S}_{n/2}}(w_{\pi})$ by conjugation. Now applying the explicit induction formula we obtain
\begin{align*}
\Theta(w_{2\pi}^+) &= \sum_{(\delta,\epsilon) \in \mathcal{B}_{\pi}} \frac{|C_{W}(w_{2\pi}^+)|}{|C_{H}(w_{2\delta}^+w_{2\epsilon}^+)|}\Delta_A(w_{2\delta}^+)\Delta_B(w_{2\epsilon}^+)\\
&\qquad +\sum_{(\delta,\epsilon) \in \mathcal{B}_{\pi}} \frac{|C_{W}(w_{2\pi}^+)|}{|C_{H}(w_{2\delta}^-w_{2\epsilon}^-)|}\Delta_A(w_{2\delta}^-)\Delta_B(w_{2\epsilon}^-)\\
&= (-1)^{n/2}2^{\ell(\pi)+1}\sum_{(\delta,\epsilon) \in \mathcal{B}_{\pi}} \frac{|C_{\mathfrak{S}_{n/2}}(w_{\pi})|}{|C_{\mathfrak{S}_{a/2}\mathfrak{S}_{b/2}}(w_{\delta}w_{\epsilon})|}[\alpha_1](w_{\delta})[\beta_1](w_{\epsilon}),\\
&= (-1)^{n/2}2^{\ell(\pi)+1}\Ind_{\mathfrak{S}_{a/2}\mathfrak{S}_{b/2}}^{\mathfrak{S}_{n/2}}([\alpha_1]\boxtimes [\beta_1])(w_{\pi}).
\end{align*}
A similar calculation for the case $w_{2\pi}^-$ shows that
\begin{equation*}
\Theta(w_{2\pi}^{\pm}) = \pm(-1)^{n/2}2^{\ell(\pi)+1}\Ind_{\mathfrak{S}_{a/2}\mathfrak{S}_{b/2}}^{\mathfrak{S}_{n/2}}([\alpha_1]\boxtimes [\beta_1])(w_{\pi}).
\end{equation*}

We are now ready to compute the multiplicity $\langle \Theta,\Delta_X\rangle_{W}$ but we first observe that
\begin{equation*}
\Theta(w_{2\pi}^+)\Delta_X(w_{2\pi}^+) = \Theta(w_{2\pi}^-)\Delta_X(w_{2\pi}^-)
\end{equation*}
for any partition $\pi \in \mathcal{P}(n/2)$, which follows immediately from the above formulas. With this we can now compute the desired multiplicity as follows
\begin{align*}
\langle \Theta,\Delta_X\rangle_{W} &= \sum_{\pi \in \mathcal{P}(n/2)} \frac{1}{|C_{W}(w_{2\pi}^+)|}\Theta(w_{2\pi}^+)\Delta_X(w_{2\pi}^+)\\
&\qquad + \sum_{\pi \in \mathcal{P}(n/2)} \frac{1}{|C_W(w_{2\pi}^-)|}\Theta(w_{2\pi}^-)\Delta_X(w_{2\pi}^-)\\
&= 4\sum_{\pi \in \mathcal{P}(n/2)} \frac{1}{|C_{\mathfrak{S}_{n/2}}(w_{\pi})|}\Ind_{\mathfrak{S}_{a/2}\mathfrak{S}_{b/2}}^{\mathfrak{S}_{n/2}}([\alpha_1]\boxtimes [\beta_1])(w_{\pi})[\gamma_1](w_{\pi})\\
&= 4\langle \Ind_{\mathfrak{S}_{a/2}\mathfrak{S}_{b/2}}^{\mathfrak{S}_{n/2}}([\alpha_1]\boxtimes [\beta_1]), [\gamma_1] \rangle_{\mathfrak{S}_{n/2}}\\
&= 4c_{\alpha_1\beta_1}^{\gamma_1}.
\end{align*}
This completes the proof.
\end{proof}

\begin{pa}
We now rephrase \cref{lem:typeD-GP-Lemm6.1.3} when $(a,b) \in \{(1,n-1),(n-1,1)\}$ to obtain the well-known type $\D$ branching rules. However first we require some notation. For any partition $\alpha = (\alpha_1,\dots,\alpha_k) \in \mathcal{P}(n)$ let $I(\alpha)$ be the set of all $d \in \{1,\dots,k\}$ such that $d = k$ or $d<k$ and $\alpha_d > \alpha_{d+1}$. For any $d \in I(\alpha)$ we denote by $\alpha^{(d)}$ the partition obtained from $\alpha$ by replacing $\alpha_d$ with $\alpha_d-1$ and removing any zero entries. Furthermore, for any bipartition $\gamma = (\gamma_1;\gamma_2) \in \mathcal{BP}(n)$ we denote by $Z(\gamma) \subseteq \mathcal{BP}(n-1)$ the set
\begin{equation*}
\{(\gamma_1^{(d)};\gamma_2),(\gamma_2;\gamma_1^{(d)}) \mid d \in I(\gamma_1)\} \cup \{(\gamma_1;\gamma_2^{(d)}),(\gamma_2^{(d)};\gamma_1) \mid d \in I(\gamma_2)\}.
\end{equation*}
The following is deduced immediately from the type $\A$ branching rules (see \cite[6.1.8]{geck-pfeiffer:2000:characters-of-finite-coxeter-groups}), Frobenius reciprocity and \cref{lem:typeD-GP-Lemm6.1.3}.
\end{pa}

\begin{cor}
Assume $n \geqslant 4$ and let $H \leqslant W_n$ be the subgroup $W_1W_{n-1}$ or $W_{n-1}W_1$ (c.f.\ \cref{pa:reflection-subgroups-WaWb}). Then for any irreducible characters $B = \overline{[\beta_1;\beta_2]}_{\pm} \in \Irr(H)$ and $X = \overline{[\gamma_1;\gamma_2]}_{\pm} \in \Irr(W_n)$ we have
\begin{equation*}
\langle\Res_H^{W_n}(X), B\rangle_H = \begin{cases}
0 &\text{if }\beta \not\in Z(\gamma)\\
1 &\text{if }\beta \in Z(\gamma)
\end{cases}
\end{equation*}
where $\beta = (\beta_1;\beta_2) \in \mathcal{BP}(n-1)$ and $\gamma = (\gamma_1;\gamma_2) \in \mathcal{BP}(n)$.
\end{cor}

\renewcommand*{\bibfont}{\small}
\begin{spacing}{0.96}
\printbibliography
\end{spacing}
\end{document}